\documentclass[11pt]{article}
\usepackage{amsfonts}
\usepackage{graphics}
\usepackage{indentfirst}
\usepackage{color}
\usepackage{cite}
\usepackage{latexsym}
\usepackage[paper=a4paper, left=2.35cm, right=2.35cm, top=2.5cm, bottom=1.8cm, headheight=5.5pt, footskip=0.8cm, footnotesep=0.8cm, centering, includefoot]{geometry}
\usepackage{amsmath}
\allowdisplaybreaks
\usepackage{amssymb}
\usepackage[dvips]{epsfig}
\usepackage{amscd}

\newtheorem{theorem}{Theorem}[section]
\newtheorem{remark}{Remark}[section]

\newtheorem{definition}{Definition}[section]
\newtheorem{lemma}[theorem]{Lemma}

\newcommand{\vp}{\varphi}
\newcommand{\n}{\rho}

\newcommand{\ti}{\tilde}
\newcommand{\mr}{\mathbb{R}}

\DeclareMathOperator{\loc}{loc}

\renewcommand{\div}{ {\rm div }  }

\newcommand{\pa}{\partial}
\renewcommand{\r}{\mathbb{R}}

\newcommand{\bt}{\begin{theorem}}
\newcommand{\bl}{\begin{lemma}}
\newcommand{\el}{\end{lemma}}
\newcommand{\et}{\end{theorem}}

\newcommand{\al}{\alpha}

\newcommand{\ve}{\varepsilon}
\newcommand{\la}{\label}

\newcommand{\ol}{\overline}

\newcommand{\bn}{\begin{eqnarray}}
\newcommand{\en}{\end{eqnarray}}
\newcommand{\bnn}{\begin{eqnarray*}}
\newcommand{\enn}{\end{eqnarray*}}

\newcommand{\bnnn}{\begin{eqnarray*}}
\newcommand{\ennn}{\end{eqnarray*}}

\newcommand{\ba}{\begin{aligned}}
\newcommand{\ea}{\end{aligned}}
\newcommand{\be}{\begin{equation}}
\newcommand{\ee}{\end{equation}}
\def\O{{\r^2 }}
\def\p{\partial}
\def\norm[#1]#2{\|#2\|_{#1}}

\def\la{\label}

\def\na{\nabla}

\makeatletter      
\@addtoreset{equation}{section}
\makeatother

\makeatletter      
\@addtoreset{equation}{section}
\makeatother       

\title{Global Existence and Large Time Asymptotic Behavior of Strong Solutions to the Cauchy Problem of 2D Density-Dependent Magnetohydrodynamic Equations with Vacuum\thanks{B. L\"u is  supported by  NNSFC Tianyuan (No. 11426131) and  Natural Science Foundation of Jiangxi Province (No. 20142BAB211006).}}

\author{ Boqiang L\"u\thanks{College of Mathematics and Information Science, Nanchang Hangkong University, Nanchang 330063, People's Republic of China ({\tt lvbq86@163.com}). }
\quad   Zhonghai Xu\thanks{College of Science, Northeast Dianli University, Jilin 132013, People's Republic of China ({\tt xuzhonghai@163.com}).
 }
 \quad Xin Zhong\thanks{Corresponding author. Institute of Applied Mathematics, AMSS,
Chinese Academy of Sciences, Beijing 100190,  People's Republic of China ({\tt xzhong1014@amss.ac.cn}).
 }
 }

\date{ }

\begin{document}
\maketitle

 \begin{abstract}
This paper concerns the Cauchy problem of the two-dimensional (2D) nonhomogeneous incompressible Magnetohydrodynamic (MHD) equations with vacuum as far field density. We establish the global existence and uniqueness of strong solutions to the 2D Cauchy problem on the whole space $\r^2$,  provided that the initial density and  the initial magnetic  decay not too  slow  at infinity.  In particular, the initial data can be arbitrarily large and the initial density can contain vacuum states and even have compact support. Furthermore,  we also obtain the large time decay rates of the gradients of velocity,   magnetic  and  pressure.
 \end{abstract}


Math Subject Classification: 35Q35; 76D03; 76W05.

\section{Introduction}

Magnetohydrodynamics is concerned with the interaction between fluid flow and magnetic field. The governing equations of nonhomogeneous incompressible MHD can be stated as follows \cite{david},
\be \la{1.1}  \begin{cases} \n_t + \div(\n u) = 0,\\
 (\n u)_t + \div(\n u\otimes u) + \nabla p = \mu\Delta u + H\cdot\na H-\frac{1}{2}\na |H|^2,\\
 H_t-\nu\Delta H+u\cdot\na H-H\cdot\na u=0,\\
 \div u=0,\,\,\, \div H=0.
\end{cases}\ee
Here,  $t\ge 0$ is time,   $x=(x_1,x_2)\in \Omega\subset \r^2$   is the spatial coordinate,  and $\rho=\rho(x,t)$, $ u =(u^1,u^2)(x,t)$, $ H =(H^1,H^2)(x,t)$, and $p=p(x,t)$ denote the density, velocity,  magnetic, and pressure of the fluid, respectively;
$\mu>0$ stands for the viscosity constant; the constant $\nu>0$ is the resistivity coefficient which is inversely proportional to the electrical conductivity constant and acts as the magnetic diffusivity of magnetic fields.

Let $\Omega=\r^2$ and we consider the Cauchy problem for \eqref{1.1} with $(\rho, u, H)$ vanishing at infinity (in some weak sense) and the initial conditions:
\begin{equation}\label{1.3'}
\rho(x,0)=\rho_0(x),\,   \n u (x,0)= \n_0u_0(x),\, H (x,0)=H_0(x), \ x\in\mathbb{R}^2,
\end{equation}
 for given initial data $\rho_0, u_0$ and $H_0$.

Magnetohydrodynamics studies the dynamics of electrically conducting fluids
and the theory of the macroscopic interaction of electrically conducting fluids with a
magnetic field. The dynamic
motion of the fluid and the magnetic field interact strongly with each other, so the
hydrodynamic and electrodynamic effects are coupled.  If this motion occurs in the absence of magnetic field, that is, $H=0$,
  the MHD system  reduces to the Navier-Stokes equations,
which have been discussed in numerous studies, please refer to  \cite{AK1973,AKM1990,CK2003,HW2014,HW2015,K1986,K1974,L1996,L2015,lvzh1,S1990,zhang,F2004,H19951,HL20130,HL20131,HLX2012,LL2014,LX2014}.  In general, due to the similarity of the second equation and the third equation in \eqref{1.1}, the study for MHD system has been along with that for Navier-Stokes one. However, the issues of well-posedness and dynamical behaviors of MHD system are rather
complicated to investigate because of the strong coupling and interplay interaction
between the fluid motion and the magnetic field.

First, let us give a short survey for the study of incompressible Navier-Stokes equations, that is,
the system \eqref{1.1} with $H=0$. In the case when $\rho_0$ is bounded away from zero,
Kazhikov \cite{K1974} established the global existence of weak solutions (see also \cite{AK1973}).
Later, Antontsev-Kazhikov-Monakhov \cite{AKM1990} gave the first result on local existence
and uniqueness of strong solutions, and then proved the unique local strong solution is
 global in two dimensions. When the initial data may contain vacuum states,
 Simon \cite{S1990} obtained the global existence of weak solutions, see also Lions \cite{L1996} for the case of density-dependent viscosity.
   Choe-Kim \cite{CK2003} proposed a compatibility condition and established
   the local existence of strong solutions.   Under some suitable smallness conditions, the global existence of strong  solutions on   bounded domains were established by Huang-Wang \cite{HW2014,HW2015}  and Zhang \cite{zhang}, respectively. Recently,  the local and global (with general  large data) existence of strong solutions to the 2D Cauchy problem with vacuum on the whole space $\O$ were established, by Liang \cite{L2015} and L\"u-Shi-Zhong \cite{lvzh1}, respectively.


Let's go back to the MHD system \eqref{1.1}. When $\n$ is a constant, which means the fluid is homogeneous, the MHD system has been extensively studied.  Duraut-Lions \cite{dlions}  constructed a class of weak solutions with finite energy and a class of local strong solutions (the local strong solution has been proved to be global in two dimensions,
 but  only local in three dimensions  except for small data), see also Sermange-Temam \cite{teman}.
  For the nonhomogeneous case, Gerbeau and Le Bris \cite{leb1}, Desjardins and
Le Bris \cite{leb2} studied the global existence of weak solutions with finite energy on 3D bounded domains  and on the torus, respectively.  In the absence of vacuum,  Abidi-Hmidi \cite{abidi1} and Abidi-Paicu \cite{abidi2} established  the local and global (with small initial data) existence of strong solutions in some Besov spaces, respectively.  In  the presence of  vacuum, under  the following  compatibility conditions,
\be\ba\label{tan}
\div u_0=\div H_0=0,\,\,\,-\Delta u_0+\na p_0-(H_0\cdot \na) H_0=\n_0^{1/2}g,~~~\mbox{in}~~\Omega,
\ea\ee
where $(p_0,~g)\in H^1\times L^2$ and $\Omega=\mathbb{R}^3$,
Chen-Tan-Wang \cite{tan} obtained the local existence of strong solutions to the 3D Cauchy problem.
 When $\Omega \subset \mathbb{R}^2$ is a bounded domain,
 Huang-Wang \cite{hwjde1} investigated  the global existence
 of strong solution with general large data when the initial density contains vacuum states and the initial data satisfy the compatibility conditions \eqref{tan}.
  Very recently,  L\"u-Xu-Zhong \cite{lvzh}  established the local existence
   of strong solution to the 2D Cauchy problem  \eqref{1.1} with vacuum as far field density.
    However, the global existence of strong solution with general large data to the 2D Cauchy problem \eqref{1.1} with  vacuum as far field density is still open. In fact, this is the main aim of this paper.

Before stating the main results, we first explain the notations and
conventions used throughout this paper.  For $R>0$, set
$$B_R  \triangleq\left.\left\{x\in\r^2\right|
\,|x|<R \right\} , \quad \int \cdot dx\triangleq\int_{\r^2}\cdot dx.$$ Moreover, for $1\le r\le \infty$ and  $k\ge 1, $  the standard Lebesgue and  Sobolev spaces are defined as follows:
   \bnn
L^r=L^r(\r^2 ),\quad
W^{k,r}  = W^{k,r}(\r^2) , \quad H^k = W^{k,2} .
 \enn

 Next, we give the definition of strong solution to \eqref{1.1} as follows:

\begin{definition}\label{def1}
If all derivatives involved in \eqref{1.1} for $(\rho, u ,p, H)$ are regular distributions, and equations \eqref{1.1} hold almost everywhere in $\mathbb{R}^2\times(0,T)$, then $(\rho, u ,p, H)$ is called a strong solution to \eqref{1.1}.
\end{definition}

Without loss of generality, we assume that the initial density $\n_0$ satisfies
\be\la{oy3.7} \int_{\r^2} \n_0dx=1,\ee  which implies that there exists a positive constant $N_0$ such that  \be\label{1.3} \int_{B_{N_0}}  \n_0  dx\ge \frac12\int\n_0dx=\frac12.\ee

Our main result can be stated as follows:
\begin{theorem}\label{thm1} In addition to \eqref{oy3.7} and \eqref{1.3}, assume that the initial data $(\rho_0, u_0, H_0)$ satisfy for any given numbers $  a>1$ and $ q>2$,
\be\ba\label{2.2}\begin{cases} \n_0\ge0,\,\rho_{0}\bar{x}^{a}\in L^{1}\cap H^{1}\cap W^{1,q},\, H_0\bar{x}^{a/2}\in L^{2},\,   \sqrt{\rho_0} u_0\in L^2,\\
\nabla u_{0} \in L^2,\ \nabla H_{0}\in L^2,\,\,\div u_0=\div H_0=0,
\end{cases}
\ea\ee
where
\begin{equation}\label{2.1}
\bar{x}\triangleq(e+|x|^2)^{1/2}\log^{2}(e+|x|^2).
\end{equation}
Then  the problem \eqref{1.1}-\eqref{1.3'} has a unique global strong solution $(\rho, u ,p, H)$ satisfying that
for any $0<T<\infty$,
\be\label{2.3}\begin{cases}
0\le\rho\in C([0,T];L^1 \cap H^1\cap W^{1,q} ),\\
\rho\bar x^a\in L^\infty( 0,T ;L^1\cap H^1\cap W^{1,q} ),\\
\sqrt{\n } u,\,\na u,\, \bar x^{-1}u,\,    \sqrt{t} \sqrt{\n}  u_t,\,    \sqrt{t} \na P,\,    \sqrt{t} \na^2  u \in L^\infty(0,T; L^2 ), \\
H,  H \bar{x}^{a/2}, \na H,  \sqrt{t}H_t,\,    \sqrt{t} \na^2 H \in L^\infty( 0,T ;L^2), \\
\na u\in  L^2(0,T;H^1)\cap  L^{(q+1)/q}(0,T; W^{1,q}), \\
\na P\in  L^2(0,T;L^2)\cap  L^{(q+1)/q}(0,T;L^q), \\
\na H\in L^2(0, T; H^1),\,\,H_t,~\na H\bar{x}^{a/2}\in L^2(0, T; L^2),\\
\sqrt{t}\na u\in L^2(0,T; W^{1,q} ),  \\
\sqrt{\n} u_t, \,\sqrt{t}\na H\bar{x}^{a/2}, \,  \sqrt{t}\na u_t ,\, \sqrt{t}\na H_t ,\,
\sqrt{t} \bar x^{-1}u_t\in L^2(\mathbb{R}^2\times(0,T)),\\
   \end{cases}\ee
and
\be\la{l1.2}
\inf\limits_{0\le t\le T}\int_{B_{N_1}}\n(x,t) dx\ge \frac14,
\ee
for some positive constant $N_1$ depending only on $\|\n_0\|_{L^1}, \|\n_0^{1/2}u_0\|_{L^2}, N_0 $, and $T$. Moreover, $(\n,u, p, H)$ has the following decay
rates,  that is,  for $t\ge 1,$
 \be \la{lv1.2}
\begin{cases}
\|\na u(\cdot,t)\|_{L^2}+\|\na H(\cdot,t)\|_{L^2}\le Ct^{-1/2},\\
\|\na^2u(\cdot,t)\|_{L^2}+\|\na p(\cdot,t)\|_{L^2}+\||H||\na H|\|_{L^2}\le Ct^{-1},
\end{cases} \ee
where 
 $C$ depends only on  $  \mu ,  \nu,  \|\n_0\|_{L^1\cap L^\infty}, \|\n_0^{1/2}u_0\|_{L^2}, \|\na u_0\|_{L^2}$, and  $\|H_0\|_{H^1}$.
\end{theorem}

\begin{remark}
\la{re1}
When there is no electromagnetic field effect, that is  $H=0$,
\eqref{1.1} turns to be the  incompressible Navier-Stokes equations,
and Theorem \ref{thm1}  is similar to the results of \cite{lvzh1}.
Roughly speaking, we generalize the results of  \cite{lvzh1} to  the  incompressible  MHD system.
Furthermore, the large time  decay rates \eqref{lv1.2} with $H=0$ are the same as those in \cite{lvzh1}, hence the magnetic field has no influence on the large time  behaviors of the velocity and the pressure.

\end{remark}



\begin{remark}\la{re5}  Our Theorem 1.1 holds for arbitrarily large initial data which is in sharp contrast to L\"u-Shi-Xu \cite{lvshixu} where the smallness conditions on the initial energy is needed in order to obtain the global existence of strong solutions to the 2D compressible MHD equations.   \end{remark}

\begin{remark}\la{re4}
Compared with \cite{hwjde1,tan}, there is no need to impose the additional compatibility condition on the initial data for the global existence of the strong solution.
\end{remark}

We now make some comments on the analysis of this paper. Note that for initial data
in the class satisfying \eqref{2.2}, the local existence
and uniqueness of strong  solutions to the Cauchy problem, (1.1)-(1.2), have been
established recently in \cite{lvzh} (see Lemma 2.1). To extend the strong solution globally in time, one
needs some global a priori estimates on strong solutions to (1.1)-(1.2) in suitable higher
norms. It should be pointed out that, on the one hand,
the crucial techniques of proofs in \cite{hwjde1,HW2014} cannot be adapted to the situation treated here,
since   their arguments only hold true for the case of bounded domains.
On the other hand,  it seems difficult to bound the  $L^q(\mathbb{R}^2)$-norm of $ u $
  just in terms of $\|\sqrt{\rho} u \|_{L^{2}(\mathbb{R}^2)}$ and $\|\nabla u \|_{L^{2}(\mathbb{R}^2)}$. 
  To this end,   we try to adapt some basic ideas used in \cite{lvzh1}, where they investigated  the global existence of strong solutions to 2D Cauchy problem of the density-dependent  Navier-Stokes equations. However, compared with \cite{lvzh1},  for the incompressible MHD equations treated here, the strong
coupling between the velocity field and the magnetic field, such as $u\cdot \na H$,   will bring out some new difficulties.

To overcome these difficulties stated  above, some new ideas are needed.
 First, we try to obtain the estimates on the $L^\infty(0,T;L^2(\r^2))$-norm  of the gradients of  velocity and  magnetic.
  On the one hand, motivated by \cite{H19951,HLX2012,LX2014},  multiplying \eqref{1.1}$_2$ by the material derivatives  $\dot u\triangleq u_t+u\cdot\na u$  instead of the usual $u_t$ (see \cite{hwjde1,HW2014}), the key point is to   control the term $\int |p||\na u|^2dx$. Motivated by \cite{lvzh1} (see also \cite{D1997}), using some facts on Hardy and $BMO$ spaces (see Lemma \ref{lem27}), we succeed in bounding the term $\int |p||\na u|^2dx$ by $ \|\na p\|_{L^2}\|\na u\|_{L^2}^2$ (see \eqref{3.8}). On the other hand,   the usual  $L^2(\O\times(0,T))$-norm  of $H_t$ cannot be directly estimated due to the  strong
coupled term  $u\cdot \na H$.   Motivated by \cite{lvshixu}, multiplying   \eqref{1.1}$_3$ by $\Delta H$   instead of the usual $H_t$ (see \cite{hwjde1}),
the  coupled term  $u\cdot \na H$ can be controlled after integration by parts (see \eqref{mm1}).
Next,  using the structure of the 2D magnetic equation (see \eqref{amss2} and \eqref{lv3.57'1}-\eqref{lv3.571}),
 we multiply \eqref{1.1}$_3$ by $H\Delta |H|^2$ and
 thus obtain some useful a priori estimates on
  $\||H||\na H|\|_{L^2}$ and   $\||H||\Delta  H|\|_{L^2}$,
 which are crucial in deriving the time-independent estimates on
 both the   $L^\infty(0,T;L^2(\r^2))$-norm of $t^{1/2} \n^{1/2}\dot u$
 and the $L^2(\r^2\times(0,T ))$-norm of $t^{1/2}\na \dot u$ (see \eqref{3.13}).
 This together with  some careful analysis on
 the spatial weighted estimates of the density (see \eqref{06.1})
  indicates the desired $L^1(0,T;L^\infty)$-bound of the gradient
  of the velocity (see \eqref{tt}) which in particular implies
   the  bound  on the $L^\infty(0,T;L^q)$-norm of the gradient of the density.
 With the a priori estimates stated above at hand, using the  similar arguments as in \cite{lvzh1,LX2014,lvshixu},
the next  step  is to bound
the higher order derivatives of the solutions $(\n, u, p, H)$.
 Finally, some useful spatial weighted  estimates on both $H$ and $\na H$ (see \eqref{gj10} and \eqref{gj10'}) are derived, such a derivation yields the estimate  on the $L^2(\r^2\times(0,T ))$-norms of both $t^{1/2}\na u_t$ and $t^{1/2}\na H_t$, and simultaneously also the bound  of the  $L^\infty(0,T; L^2(\r^2))$-norm of  $t^{1/2}\na^2 H$, see Lemma \ref{lem4.5v} and its proof.

The rest of this paper is organized as follows. In Section \ref{sec2}, we collect some elementary facts and inequalities that will be used later. Section \ref{sec3} is devoted to the a priori estimates. Finally,  Theorem \ref{thm1} is proved in Section \ref{sec4}.

\section{Preliminaries}\label{sec2}

In this section, we will recall some  known facts and elementary
inequalities which will be used frequently later.
We start with the local existence theorem of strong solutions whose proof can be found in \cite{lvzh}.
\begin{lemma}\label{lem21}
Assume that $(\rho_0, u_0, H_0)$ satisfies \eqref{2.2}. Then there exists a small time $T>0$ and a unique strong solution $(\rho, u, p, H )$ to the problem \eqref{1.1}-\eqref{1.3'} in $\mathbb{R}^{2}\times(0,T)$ satisfying \eqref{2.3} and \eqref{l1.2}.
\end{lemma}

Next, the following well-known Gagliardo-Nirenberg inequality (see \cite{nir})
  will be used later.

\begin{lemma}
[Gagliardo-Nirenberg]\la{l1} For  $s\in [2,\infty),q\in(1,\infty), $ and
$ r\in  (2,\infty),$ there exists some generic
 constant
$C>0$ which may depend  on $s,q, $ and $r$ such that for   $f\in H^1({\O }) $
and $g\in  L^q(\O )\cap D^{1,r}(\O), $    we have \be
\la{g1}\|f\|_{L^s(\O)}^s\le C \|f\|_{L^2(\O)}^{2}\|\na
f\|_{L^2(\O)}^{s-2} ,\ee  \be
\la{g2}\|g\|_{C\left(\ol{\O }\right)} \le C
\|g\|_{L^q(\O)}^{q(r-2)/(2r+q(r-2))}\|\na g\|_{L^r(\O)}^{2r/(2r+q(r-2))} .
\ee
\end{lemma}

The following weighted $L^m$ bounds for elements in $\tilde{D}^{1,2}(\O)\triangleq\{ v \in H_{\loc}^{1}(\mathbb{R}^2)|\nabla v \in L^{2}(\mathbb{R}^2)\}$ can be found in \cite[Theorem B.1]{L1996}.
\begin{lemma} \la{1leo}
For $m\in [2,\infty)$ and $\theta\in (1+m/2,\infty),$ there exists a positive constant $C$ such that we have for all $v\in  \tilde{D}^{1,2}(\O),$ \be\la{3h} \left(\int_{\O} \frac{|v|^m}{e+|x|^2}(\log (e+|x|^2))^{-\theta}dx  \right)^{1/m}\le C\|v\|_{L^2(B_1)}+C\|\na v\|_{L^2(\O) }.\ee
\end{lemma}

The combination of Lemma \ref{1leo} and the Poincar\'e inequality yields
 the following useful results on  weighted bounds, whose proof can be found  in \cite[Lemma 2.4]{LX2014}.

\begin{lemma}\label{lem26} Let  $\bar x$  be as in \eqref{2.1}. Assume that $\n \in L^1(\O)\cap L^\infty(\O)$ is a non-negative function such that
\be \la{2.i2}   \int_{B_{N_1} }\n dx\ge M_1,  \quad \|\n\|_{L^1(\O)\cap L^\infty(\O)}\le M_2,  \ee
for positive constants $   M_1, M_2, $ and $ N_1\ge 1$  with $B_{N_1}\subset\O.$ Then for $\ve> 0$ and $\eta>0,$ there is a positive constant $C$ depending only on   $\ve,\eta, M_1,M_2, $ and $N_1$ such that      every $v\in \ti D^{1,2}(\O)  $ satisfies\be\la{3.i2}\ba \|v\bar x^{-\eta}\|_{L^{(2+\ve)/\ti\eta}(\O)} &\le C \|\n^{1/2}v\|_{L^2(\O)}+C \|\na v\|_{L^2(\O)} ,\ea\ee with $
\ti\eta=\min\{1,\eta\}.$
\end{lemma}

Finally, let $\mathcal{H}^{1}(\mathbb{R}^2)$ and $BMO(\mathbb{R}^2)$ stand for the usual Hardy and
$BMO$ spaces (see \cite[Chapter IV]{S1993}). Then the following well-known facts play a key role in the proof of Lemma \ref{lem3.3} in the next section.
\begin{lemma}\label{lem27}
(a) There is a positive constant $C$ such that
\begin{equation*}
\|E\cdot B\|_{\mathcal{H}^{1}(\O)}
\leq C\|E\|_{L^{2}(\O)}\|B\|_{L^{2}(\O)},
\end{equation*} for all $E\in L^{2}(\mathbb{R}^2)$ and $ B\in L^{2}(\mathbb{R}^2)$ satisfying
\begin{equation*}
\div E=0,\ \nabla^{\bot}\cdot B=0\ \ \text{in}\ \ \mathcal{D}'(\mathbb{R}^2).
\end{equation*}
(b) There is a positive constant $C$ such that
\begin{equation}\label{lem1}
\| v \|_{BMO(\O)}\leq C\|\nabla v \|_{L^{2}(\O)},
\end{equation} for all $  v \in D^1(\mathbb{R}^2)$.
\end{lemma}
{\it Proof.}
(a) For the detailed proof, please see \cite[Theorem II.1]{CLMS1993}.

(b) It follows  from the Poincar{\'e} inequality that for any ball $B\subset\mathbb{R}^2$
\begin{equation*}
\frac{1}{|B|}\int_{B}\left| v(x) - \frac{1}{|B|}\int_Bv(y)dy\right|dx\leq C\left(\int_{B}|\nabla v |^2dx\right)^{1/2},
\end{equation*}
which directly gives \eqref{lem1}.\hfill $\Box$

\section{A Priori Estimates}\label{sec3}
In this section, we will establish some necessary a priori bounds for strong solutions $(\rho, u ,p, H)$ to the Cauchy problem \eqref{1.1}-\eqref{1.3'} to extend the local strong solutions guaranteed by Lemma \ref{lem21}. Thus, let $T>0$ be a fixed time and $(\rho,u,p,H)$ be the strong solution to \eqref{1.1}-\eqref{1.3'} on $\mathbb{R}^{2}\times(0,T]$ with initial data $(\rho_0, u_0, H_0)$ satisfying \eqref{oy3.7}-\eqref{2.2}.

In what follows, we will use the convention that $C$ denotes a generic positive constant depending on $\mu$, $\nu$, $a$, and the initial data,  and  use $C(\al)$ to emphasize that $C$ depends on $\al.$

\subsection{Lower Order Estimates}
First, since $\div u=0,$ we have the following estimate on the $L^\infty(0,T;L^r)$-norm of the density.
\begin{lemma}[\cite{L1996}]\label{lem3.1}
There exists a positive constant $C$ depending only on $\|\n_0\|_{L^1\cap L^\infty}$ such that
\begin{equation}\label{3.1}
\sup_{t\in[0,T]}\|\rho\|_{L^{1}\cap L^\infty}\leq C.
\end{equation}
\end{lemma}

The following lemma concerns the time-independent estimates on the $L^\infty(0,T;L^2)$-norm of the gradients of the velocity and the magnetic.
\begin{lemma}\label{lem3.3}
There exists a positive constant $C$ depending only on $  \mu ,$  $\nu, $ $ \|\n_0\|_{L^\infty},$ $  \|\n_0^{1/2}u_0\|_{L^2},$ $\|\na u_0\|_{L^2}$,  and  $\|H_0\|_{H^1}$ such that
\begin{align}\label{lv3.5}
&\sup_{t\in[0,T]}\left( \|\nabla u\|_{L^2}^2+\|\na H\|_{L^2}^2+\|H\|_{L^4}^4\right) \notag \\
&\quad + \int_0^{T}
\left(\|\n^{1/2}\dot{u} \|_{L^2}^2 +\|\triangle H\|_{L^2}^2+\||\na H| |H|\|_{L^2}^2\right) dt\leq C.
\end{align}
Here $\dot{v}\triangleq\partial_{t}v+u\cdot\nabla v$. Furthermore, one has
\begin{align}\label{lv3.5'}
&\sup_{t\in[0,T]}t\left(\|\nabla u\|_{L^2}^2+\|\na H\|_{L^2}^2+\|H\|_{L^4}^4\right) \notag \\
&\quad+ \int_0^{T} t
\left(\|\n^{1/2}\dot{u}\|_{L^2}^2+\|\triangle H\|_{L^2}^2+\||\na H| |H|\|_{L^2}^2\right) dt\leq C.
\end{align}
\end{lemma}

{\it Proof.}
First, applying standard energy estimate to \eqref{1.1} gives
\be\ba\la{lvgj1}
\sup_{t\in[0,T]} \left( \|\rho^{1/2}u\|^2_{L^{2}}+\|  H\|_{L^2}^2\right)
+\int_{0}^{T}\left( \|\na u\|_{L^2}^2+\|\na H\|_{L^2}^2 \right)dt  \le C.
\ea\ee

Next, multiplying \eqref{1.1}$_2$ by $\dot{u}$ and integrating the resulting equality over $\mathbb{R}^2$ lead  to
\begin{align}\label{3.6}
\int\rho|\dot{ u }|^{2}dx&=\int \mu\Delta u \cdot\dot{ u }dx-\int\nabla p\cdot\dot{ u } dx-\frac{1}{2}\int \na |H|^2\cdot\dot udx+\int H\cdot \na H\cdot\dot u dx \notag \\
&\triangleq I_{1}+I_{2}+I_{3}+I_{4}.
\end{align}
It follows from integration by parts and Gagliardo-Nirenberg inequality that
\begin{align}\label{3.7}
I_{1} & =\int\mu\Delta u \cdot( u_t  + u \cdot\nabla u )dx \notag \\
& = -\frac{\mu}{2}\frac{d}{dt}\|\nabla u \|_{L^2}^{2}
-\mu\int\partial_{i}u^{j}\partial_{i}(u^{k}\partial_{k}u^{j})dx  \notag \\
& \leq-\frac{\mu}{2}\frac{d}{dt}\|\nabla u \|_{L^2}^{2}+C\|\nabla u \|_{L^3}^{3} \notag \\
& \leq-\frac{\mu}{2}\frac{d}{dt}\|\nabla u \|_{L^2}^{2}
+C\|\nabla u \|_{L^2}^{2}\|\nabla^{2} u \|_{L^2}.
\end{align}
Integration by parts together with \eqref{1.1}$_4$ gives
\begin{align}\label{mhd3.7}
I_{2}& =-\int\nabla p\cdot( u_t  + u \cdot\nabla u )dx \notag \\
& = \int p\partial_{j}u^{i}\partial_{i}u^{j}dx \notag \\
& \leq C\|p\|_{BMO}\|\partial_{j}u^{i}\partial_{i}u^{j}\|_{\mathcal{H}^1},
\end{align}
where one has used the duality of $\mathcal{H}^1$ and $BMO$ (see \cite[Chapter IV]{S1993}) in the last inequality. Since $\div(\partial_j u )=\partial_j\div u =0$ and $\nabla^{\bot}\cdot(\nabla u^{j})=0$,   Lemma \ref{lem27} yields
\begin{equation}\label{3.8}\ba
|I_{2}|=\left|\int p\partial_{j}u^{i}\partial_{i}u^{j}dx \right|\leq C\|\nabla p\|_{L^2}\|\nabla u \|_{L^2}^{2}.\ea
\end{equation}
For the term $I_3$, integration by parts together with \eqref{1.1}$_4$ and \eqref{g1} leads to
\begin{align}\label{lib1}
I_3
& = \frac{1}{2}\int |H|^2 \p_iu^j\p_ju^i dx \notag \\
& \leq C\|H\|_{L^6}^6+C\|\na u\|_{L^3}^3 \notag \\
& \leq  C\|H\|_{L^2}^2\|\na H\|_{L^2}^4+C\|\nabla u \|_{L^2}^{2}\|\nabla^{2} u \|_{L^2}.
\end{align}
Using \eqref{1.1}$_3$ and \eqref{1.1}$_4$, one deduces from integration by parts and \eqref{g1} that
\begin{align}\la{lib2}
I_4
& = \int H\cdot\na H\cdot u_tdx+\int H\cdot\na H\cdot (u\cdot \na u)dx
\notag \\
& = -\frac{d}{dt}\int H\cdot\na u \cdot  Hdx+\int H_t\cdot\na u\cdot Hdx+\int H\cdot\na u\cdot  H_tdx \notag \\
& \quad- \int H^i\p_i u^j\p_j u^k H^kdx- \int H^i u^j\p_i\p_j u^k H^kdx
\notag \\
& =-\frac{d}{dt}\int H\cdot\na u \cdot  Hdx+\int (\nu\Delta H-u\cdot\na H+H\cdot\na u)\cdot\na u\cdot  Hdx \notag \\
& \quad+\int H\cdot\na u\cdot (\nu\Delta H-u\cdot\na H+H\cdot\na u)dx- \int H^i\p_i u^j\p_j u^k H^kdx \notag \\
& \quad+\int u^j\p_jH^i \p_i u^k H^kdx +\int H^i \p_i u^k u^j\p_jH^kdx
\notag \\
& \le -\frac{d}{dt}\int H\cdot\na u \cdot  Hdx+\frac\nu2\|\Delta H\|_{L^2}^2+C\|  H\|_{L^2}^2\|\na H\|_{L^2}^4 +C\|\nabla u \|_{L^2}^{2}\|\nabla^{2} u \|_{L^2}.
\end{align}
Hence, inserting \eqref{3.7} and \eqref{3.8}-\eqref{lib2} into \eqref{3.6} indicates that
\begin{align}\label{3.9}
& \frac{d}{dt}\left(\frac{\mu}{2}\|\nabla u \|_{L^2}^{2}+\int H\cdot\na u \cdot  Hdx\right)+\|\sqrt{\rho}\dot{ u }\|_{L^2}^{2} \notag \\
& \leq \frac\nu2\|\Delta H\|_{L^2}^2+C\|  H\|_{L^2}^2\|\na H\|_{L^2}^4
+C\left(\|\nabla^2 u \|_{L^2}+\|\nabla p\|_{L^2}\right)\|\nabla u \|_{L^2}^{2}.
\end{align}

Now, multiplying \eqref{1.1}$_3$ by $ \triangle H$ and integrating the resulting equality by parts over $\mr^2$,
it follows from H\"older's and Gagliardo-Nirenberg inequalities that
\begin{align}\la{mm1}
&\quad\frac{d}{{ d}t}\int|\na H|^2 dx+2\nu\int|\triangle H|^2dx \notag \\
&\quad\leq C\int |\na u||\na H|^2dx+ C\int |\na u|| H| |\Delta H|dx \notag \\
&\quad\leq C\|\na u\|_{L^3} \|\na H\|_{L^2}^{4/3} \|\Delta H\|_{L^2}^{2/3}+C\|\na u\|_{L^3}\| H\|_{L^6}  \|\na^2 H\|_{L^2}\notag \\
&\quad\leq C \|\na u\|_{L^2}^2\|\na^2 u\|_{L^2}+C (1+\|H\|_{L^2}^2)\|\na H\|_{L^2}^4+\frac\nu2 \|\Delta H\|_{L^2}^{2},
\end{align}
which together with \eqref{3.9} and \eqref{lvgj1} gives
\begin{align}\label{lv3.9}
& \frac{d}{dt}\left(\frac{\mu}{2}\|\nabla u \|_{L^2}^{2}+\|\nabla H\|_{L^2}^{2}+\int H\cdot\na u \cdot Hdx\right)+\|\sqrt{\rho}\dot{ u }\|_{L^2}^{2}+\nu\|\Delta H\|_{L^2}^{2} \notag \\
& \leq C\|\na H\|_{L^2}^4
+C\left(\|\nabla^2 u\|_{L^2}+\|\nabla p\|_{L^2}\right)\|\nabla u \|_{L^2}^{2}.
\end{align}

On the other hand, since $(\n,u,p,H)$ satisfies the following Stokes system
\be\la{stokes1}
\begin{cases}
-\mu\Delta u + \nabla p = -\n \dot u + H\cdot\na H-\frac{1}{2}\na |H|^2,\,\,\,\,&x\in \mathbb{R}^2,\\
\div u=0,   \,\,\,&x\in  \mathbb{R}^2,\\
u(x)\rightarrow0,\,\,\,\,&|x|\rightarrow\infty,
\end{cases}
\ee
applying the standard $L^r$-estimate to \eqref{stokes1} (see \cite{T2001}) yields that for any $r>1$,
\be\ba\label{stokes2}
\|\na^2 u \|_{L^r}+\|\nabla p\|_{L^r}& \le C\|\n \dot u \|_{L^r} + C\||H||\na H|\|_{L^r}\leq C\|\sqrt{\rho}\dot{ u }\|_{L^r}+ C\||H||\na H|\|_{L^r},
\ea\ee
where in the last inequality one has used \eqref{3.1}.

Thus, it follows from \eqref{lv3.9} and \eqref{stokes2} that
\begin{align}\label{3.11}
&\frac{d}{dt}B(t)
+\|\sqrt{\rho}\dot{ u }\|_{L^2}^{2}+\nu\|\Delta H\|_{L^2}^{2} \notag \\
&\leq  C \|\na H\|_{L^2}^4+C\|\nabla u \|_{L^2}^{4}
 +\ve\|\sqrt{\rho}\dot{ u }\|_{L^2}^{2}+ \ve\||H||\na H|\|_{L^2}^2,
\end{align}
where
\be\ba \label{lv3.11}
B(t)\triangleq \frac{\mu}{2}\|\nabla u \|_{L^2}^{2}+\|\nabla H\|_{L^2}^{2}+\int H\cdot\na u \cdot  Hdx
\ea\ee
satisfies
\be\ba\label{lllv3.11}
&\frac{\mu}{4}\|\nabla u \|_{L^2}^{2}+\|\nabla H\|_{L^2}^{2} -C_1\|H\|_{L^4}^{4} \le B(t)
\le C \|\nabla u \|_{L^2}^{2}+C\|\nabla H\|_{L^2}^{2}
\ea\ee
owing to \eqref{g1}, \eqref{lvgj1} and the following estimate
\be\ba\label{llv3.11}
\int |H\cdot\na u \cdot  H|dx\le  \frac{\mu}{4}\|\nabla u \|_{L^2}^{2}+C_1\|H\|_{L^4}^{4}.
\ea\ee

Next, multiplying $\eqref{1.1}_{3}$ by $ H|H|^2$ and integrating the resulting equality by parts over $\mr^2$ lead to
\begin{align}\la{lv3.7}
&\frac{1}{4}\left(\||H|^2\|^2_{L^2}\right)_t +\nu\||\na H| |H|\|_{L^2}^2+\frac{\nu}{2}\|\na |H|^2 \|_{L^2}^2 \notag \\
& \leq C\|\na u\|_{L^2} \||H|^2\|_{L^4}^2 \notag \\
&\leq C\|\na u\|_{L^2} \||H|^2\|_{L^2} \|\na |H|^2\|_{L^2} \notag \\
&\leq \frac{\nu}{4} \|\na |H|^2\|_{L^2}^2+ C\|\na u\|_{L^2}^4+ C\|\na H\|_{L^2}^4
\end{align}
due to Gagliardo-Nirenberg inequality and \eqref{lvgj1}.

Now, adding \eqref{lv3.7} multiplied by $4(C_1+1)$ to \eqref{3.11} and choosing $\ve$ suitably small, we obtain after using \eqref{lllv3.11} that
\begin{align}\label{lv3.12}
&\frac{d}{dt}\left(B(t)+ (C_1+1)\|H\|_{L^4}^{4}\right)
+\frac{1}{2}\|\sqrt{\rho}\dot{ u }\|_{L^2}^{2}+\frac{\nu}{2}\|\Delta H\|_{L^2}^{2}+2\nu\||H||\na H|\|_{L^2}^2 \notag \\
& \leq  C \left(\|\na H\|_{L^2}^2+\|\nabla u \|_{L^2}^{2}\right)\left(B(t)+ (C_1+1)\|H\|_{L^4}^{4}\right).
\end{align}
This combined with \eqref{lvgj1}, \eqref{lllv3.11}, and Gronwall's inequality gives \eqref{lv3.5}.

Finally, applying Gronwall's inequality to \eqref{lv3.12} multiplied by $t$, together with \eqref{lvgj1} and \eqref{lllv3.11} yields \eqref{lv3.5'} and
finishes the proof of Lemma \ref{lem3.3}.  \hfill $\Box$

\begin{lemma}\label{lem3.4}
There exists a positive constant $C$ depending only on
$ \mu , \nu, $ $\|\n_0\|_{L^1\cap L^\infty}, $ $\|\n_0^{1/2}u_0\|_{L^2},$
$ \|\na u_0\|_{L^2}$, and  $\|H_0\|_{H^1}$  such that for $i=1,2,$
\begin{equation}\label{3.13}
\sup_{t\in[0,T]}t^i\left(\|\sqrt{\rho}\dot{ u }\|_{L^2}^2+\||H||\na H| \|_{L^2}^2\right)
+\int_{0}^{T}t^i\left(\|\nabla\dot{ u }\|_{L^2}^{2} +\||\Delta H| |H|\|_{L^2}^2\right)dt\leq C,
\end{equation}
 and
\begin{equation}\label{i3.13}
\sup_{t\in[0,T]}t^i\left(\|\na^2 u\|_{L^2}^2+\|\na p\|_{L^2}^2\right) \leq C.
\end{equation}
\end{lemma}

{\it Proof.}
Motivated by \cite{H19951,LX2014,HLX2012}, operating $\partial_{t}+ u \cdot\nabla$ to \eqref{1.1}$_2^j$, one gets by some simple calculations   that
\begin{align}\label{jia4}
& \partial_{t}(\rho\dot{u}^j)+\div(\rho u \dot{u}^j)-\mu\Delta\dot{u}^j
\notag \\
& = -\mu\partial_{i}(\partial_{i} u \cdot\nabla u^{j})
-\mu\div({\partial_{i} u \partial_{i}u^j})
-\partial_{j}\partial_{t}p-( u \cdot\nabla)\partial_{j}p \notag \\
& \quad-\frac12\p_t(\p_j|H|^2)-\frac12u\cdot\na(\p_j|H|^2) +\p_t(H\cdot\na H^j)+u\cdot\na(H\cdot\na H^j),
\end{align}
which multiplied by $\dot{u}^j$, together with integration by parts and \eqref{1.1}$_4$, leads to
\begin{align}\label{3.15}
\frac{1}{2}\frac{d}{dt}\int\rho|\dot{ u }|^{2}dx+\mu\int|\nabla\dot{ u }|^{2}dx
& = -\int\mu\partial_{i}(\partial_{i} u \cdot\nabla u^{j})\dot{u}^jdx
-\int\mu\div({\partial_{i} u \partial_{i}u^j})\dot{u}^jdx \notag \\
& \quad-\int\left(\dot{u}^j\partial_{t}\partial_{j}p+\dot{u}^j( u \cdot\nabla)\partial_{j}p\right)dx \notag \\
& \quad-\frac12\int\dot{u}^j\left(\p_t(\p_j|H|^2)
+u\cdot\na(\p_j|H|^2)\right)dx
\notag \\
& \quad+\int\dot{u}^j\left(\p_t(H\cdot\na H^j)+u\cdot\na(H\cdot\na H^j)\right)dx \notag \\
&\triangleq\sum_{i=1}^{5}J_i.
\end{align}

We estimate each term on the right-hand side of \eqref{3.15} as follows.

First, by the same arguments as in \cite[Lemma 3.3]{lvzh1}, one has
\be\ba\label{i46}
\sum_{i=1}^{3}J_i\leq & \frac{d}{dt}\int p\partial_{j}u^{i}\partial_{i}u^{j}dx
+C \left(\|p\|_{L^4}^{4}+\|\nabla u \|_{L^4}^{4}\right)
+ \frac{\mu}{6}\|\nabla\dot{ u }\|_{L^2}^{2}.
\ea\ee
Next, it follows from integration by parts, \eqref{1.1}$_3$, \eqref{1.1}$_4$, and \eqref{g1} that
\begin{align}\label{lv3.44}
J_4
& = \int\pa_j\dot{u}^j H\cdot H_tdx -\frac{1}{2}\int\dot{u}^j u^i\pa_i\pa_j|H|^2 dx \notag \\
& = \int\pa_j\dot{u}^j H\cdot (H\cdot \na u+\nu \Delta H-u\cdot\na H)dx +\frac{1}{2}\int u^i\pa_i\dot{u}^j \pa_j|H|^2 dx \notag \\
& = \int\div\dot{u} H\cdot (H\cdot \na u+\nu \Delta H)dx -\frac12\int \pa_j\dot{u}^ju^i\p_i |H|^2dx+\frac{1}{2}\int u^i\pa_i\dot{u}^j \pa_j|H|^2 dx \notag \\
&\le C\int|\na u|^3 |H|^2dx +C\int|\na u|^2 |\Delta H||H| dx -\frac{1}{2}\int \pa_ju^i\pa_i\dot{u}^j |H|^2 dx \notag \\
&\le C\||H|^2\|_{L^4}^4+C\|\na u\|_{L^4}^4+C\||H||\Delta H|\|_{L^2}^2+\frac{\mu}{6}\|\na \dot u\|_{L^2}^2.
\end{align}
Similar to \eqref{lv3.44}, we also have
\begin{align}\label{lv3.45}
J_5 &=-\nu\int \Delta H\cdot\na \dot u\cdot Hdx-\nu\int H\cdot \na \dot u\cdot\Delta Hdx-\int H\cdot \na \dot u\cdot(H\cdot \na u)dx \notag \\
&\leq C\||H|^2\|_{L^4}^4+C\|\na u\|_{L^4}^4+C\||H||\Delta H|\|_{L^2}^2+\frac{\mu}{6}\|\na \dot u\|_{L^2}^2.
\end{align}
Substituting \eqref{i46}-\eqref{lv3.45} into \eqref{3.15} gives
\begin{align}\label{i47}
\frac{1}{2}\frac{d}{dt}\int\rho|\dot{ u }|^{2}dx+\frac{\mu}{2}\int|\nabla\dot{ u }|^{2}dx&\le \frac{d}{dt}\int p\partial_{j}u^{i}\partial_{i}u^{j}dx
+C\|p\|_{L^4}^{4}+C\|\nabla u \|_{L^4}^{4} \notag \\
&\quad+C\||H|^2\|_{L^4}^4 +C_2\||H||\Delta H|\|_{L^2}^2.
\end{align}

Next, as in \cite{lvhuang,lvshixu}, for $ a_1,a_2\in\{-1,0, 1\},$ denote
\be\ba\label{amss1}  \tilde{H}(a_1,a_2)=a_1H^1+a_2H^2,\quad\tilde{u}(a_1,a_2)=a_1u^1+a_2u^2,
\ea\ee
it is easy to deduce from \eqref{1.1}$_3$ that
\be\ba\label{amss2}
\tilde{H}_t-\nu\Delta \tilde{H}=H\cdot\na\tilde{u}-u\cdot \na\tilde{H}.
\ea\ee
Integrating \eqref{amss2} multiplied by $4\nu^{-1}\tilde{H} \triangle |\tilde{H}|^2$  by parts over $\mr^2$  leads to
\begin{align}\label{lv3.51}
& \nu^{-1}\left(\|\na |\tilde{H}|^2\|^2_{L^2}\right)_t+{2}\|\Delta |\tilde{H}|^2\|^2_{L^2} v \notag \\
& =4\int|\na \tilde{H}|^2 \Delta |\tilde{H}|^2dx-4\nu^{-1}\int H\cdot\na \tilde{u}\cdot\tilde{H}\Delta |\tilde{H}|^2dx \notag \\
& \quad + 2\nu^{-1}\int u\cdot \na |\tilde{H}|^2 \Delta |\tilde{H}|^2dx
\notag \\
& \le  C\|\na u\|^4_{L^4}+C \|\na H\|^4_{L^4}+C  \||H|^2\|^4_{L^4}+ \|\Delta |\tilde{H}|^2\|^2_{L^2}.
\end{align}
Noticing that
\begin{align}\label{lv3.571}
\||\Delta H||H| \|^2_{L^2}\le
& C \|\na H\|^4_{L^4}+  \|\Delta |\tilde{H}(1,0)|^2\|^2_{L^2}+ \|\Delta |\tilde{H}(0,1)|^2\|^2_{L^2} \notag \\
& + \|\Delta |\tilde{H}(1,1)|^2\|^2_{L^2}+  \|\Delta |\tilde{H}(1,-1)|^2\|^2_{L^2},
\end{align}
 and
\be\ba\label{lv3.57'1}
\||\na H||H| \|^2_{L^2}\le & G(t) \le C\||\na H||H| \|^2_{L^2} \ea\ee
with
$$
G(t) \triangleq \|\na|\tilde{H}(1,0)|^2\|^2_{L^2}+ \|\na|\tilde{H}(0,1)|^2\|^2_{L^2}  + \|\na|\tilde{H}(1,1)|^2\|^2_{L^2}+ \|\na|\tilde{H}(1,-1)|^2 \|^2_{L^2},$$
it thus follows from \eqref{lv3.51} multiplied by $(C_2+1)$ that
\begin{align}\label{lv3.51'}
& \frac{d}{dt}\left(\nu^{-1}(C_2+1)G(t)\right)+(C_2+1)\||\Delta H||H| \|^2_{L^2} \notag \\
&\le C\|\na u\|^4_{L^4}+C \|\na H\|^4_{L^4}+C\||H|^2\|^4_{L^4}.
\end{align}
This combined with \eqref{i47} yields that
\be\label{hd1}\ba
& \frac{d}{dt} F(t)
+\frac{\mu}{2}\|\na\dot u\|^2_{L^2}+ \||\Delta H||H| \|^2_{L^2}\le C\|p\|_{L^4}^{4}+C\|\nabla u \|_{L^4}^{4}+C\|\na H\|^4_{L^4}+C  \||H|^2\|^4_{L^4},
\ea\ee
where
\be\label{i98}\ba F(t) \triangleq  \frac{1}{2}\|\sqrt{\rho}\dot{ u }\|_{L^2}^{2}+\nu^{-1}(C_2+1)G(t)-\int p\partial_{j}u^{i}\partial_{i}u^{j}dx
\ea\ee
satisfies
 \be\label{i96}\ba & \frac{1}{4}\|\sqrt{\rho}\dot{ u }\|_{L^2}^{2}+\frac{\nu^{-1}(C_2+1)}{2}G(t)-C\|\nabla u \|_{L^2}^{4} \le F(t) \le   C\|\sqrt{\rho}\dot{ u }\|_{L^2}^{2}+CG(t)+C \|\nabla u \|_{L^2}^{4}
\ea\ee
owing to the following estimate
\begin{align}\label{i97}
\left|\int p\partial_{j}u^{i}\partial_{i}u^{j}dx \right|
&\le C(\|\sqrt{\n}\dot u\|_{L^2}+\||H||\na H|\|_{L^2})\|\nabla u \|_{L^2}^{2} \notag \\
&\le  \frac{1}{2}\|\sqrt{\n}\dot u\|_{L^2}^2+\frac{\nu^{-1}(C_2+1)}{2}G(t) +C\|\nabla u\|_{L^2}^{4},
\end{align}
which is deduced from   \eqref{3.8}, \eqref{stokes2}, \eqref{lv3.57'1} and Young's inequality.

Now, we shall estimate the terms on the right-hand side of \eqref{hd1}.
On the one hand, it follows from Sobolev's inequality, \eqref{stokes2}, \eqref{g1},  \eqref{3.1},   and \eqref{lvgj1} that
\begin{align}\label{11}
\|p\|_{L^4}^{4}+\|\nabla u \|_{L^4}^{4}&\leq C\|\nabla p\|_{L^{4/3}}^{4}+C\|\nabla u\|_{L^{4/3}}^{4} \notag \\
&\leq C\|\rho\dot{ u }\|_{L^{4/3}}^{4}+C\||H||\na H|\|_{L^{4/3}}^{4} \notag \\
&\leq C\|\rho\|_{L^2}^{2}\|\sqrt{\rho}\dot{ u }\|_{L^2}^{4}+ C\|H\|_{L^2}^{4}\|\na H\|_{L^4}^{4} \notag \\
&\leq C\|\sqrt{\rho}\dot{ u }\|_{L^2}^{4}+C\|\na H\|^2_{L^2}\|\na^2 H\|^2_{L^2}.
\end{align}
On the other hand, it holds from \eqref{g1} and \eqref{lvgj1} that
\be\ba\label{hd2}
\|\na H\|^4_{L^4}+ \||H|^2\|^4_{L^4}\le C\|\na H\|^2_{L^2}\|\na^2 H\|^2_{L^2}+ C\|\na H\|^2_{L^2}\||H||\na H|\|^2_{L^2}.
\ea\ee
Thus, putting \eqref{11}-\eqref{hd2} into \eqref{hd1}, together with \eqref{lv3.57'1} and \eqref{i96}, one has
\begin{align}\label{i48}
& \frac{d}{dt}F(t) +\frac{\mu}{2}\|\na\dot u\|^2_{L^2}+ \||\Delta H||H| \|^2_{L^2} \notag \\
&\le C\left(\|\sqrt{\rho}\dot{ u }\|_{L^2}^{2} + \|\na H\|^2_{L^2}\right)\left(F(t)+\|\nabla u \|_{L^2}^{4} \right)+C\|\na H\|^2_{L^2}\|\na^2 H\|^2_{L^2}.
\end{align}
Then, applying Gronwall's inequality to \eqref{i48} multiplied by $ t^i~(i=1,2)$, it follows from \eqref{lv3.57'1}, \eqref{i96}, \eqref{lv3.5},  \eqref{lv3.5'}, and \eqref{lvgj1} that
\begin{align}\label{331}
&\sup_{t\in[0,T]}\left(t^iF(t)\right)
+\int_{0}^{T}t^i\|\nabla\dot{ u }\|_{L^2}^{2}dt+\int_{0}^{T}t^i\||\Delta H||H| \|^2_{L^2}dt \notag \\
&\leq  C\int_0^T t^{i-1} F(t)ds+C\int_0^T t^i\|\na H\|^2_{L^2}\|\na^2 H\|^2_{L^2}dt \notag \\
&\quad  +C\int_0^T \left(\|\sqrt{\rho}\dot{ u }\|_{L^2}^{2} + \|\na H\|^2_{L^2}\right) t^i\|\nabla u \|_{L^2}^{4}dt \notag \\
&\leq  C\int_0^T t^{i-1} \left(\|\sqrt{\rho}\dot{ u }\|_{L^2}^{2}+\||H||\na H|\|_{L^2}^2\right)ds+C\sup_{t\in[0,T]}(t^{i-1} \|\na u\|^2_{L^2})\int_0^T  \|\na u\|^2_{L^2}dt \notag \\
&\quad+C\sup_{t\in[0,T]}\left(t^{i-1} \|\na H\|^2_{L^2}\right)\int_0^T t\|\na^2 H\|^2_{L^2}dt \notag \\
&\quad  +C\sup_{t\in[0,T]}\left(t^i \|\na u\|^4_{L^2}\right) \int_0^T \left(\|\sqrt{\rho}\dot{ u }\|_{L^2}^{2} + \|\na H\|^2_{L^2}\right) dt \notag \\
&\le C.
\end{align}
This together with \eqref{lv3.57'1}, \eqref{i96}, and \eqref{lv3.5'} yields  the desired result \eqref{3.13},
which combined with \eqref{stokes2} implies \eqref{i3.13}. The proof of Lemma \ref{lem3.4} is completed.   \hfill $\Box$

\subsection{Higher order estimates}
The following spatial weighted estimate on the density plays an important role in deriving the bounds on the higher order derivatives of the solutions $(\n,u,p,H)$.

\begin{lemma}\label{lem03.6}
There exists a positive constant $C$ depending on  $T$  such that
\begin{equation}\label{06.1}
\sup_{t\in[0,T]}\|\rho\bar{x}^{a}\|_{L^{1}}\leq C(T).
\end{equation}
\end{lemma}

{\it Proof.}
First, for $N>1$, let $\vp_N\in C^\infty_0(B_N)$  satisfy
 \be\ba \la{vp1}0\le \vp_N \le 1, \quad  \vp_N(x)=\begin{cases} 1,~~~~ |x|\le N/2,\\
0,~~~~ |x|\ge N,\end{cases}   \quad |\na  \vp_N|\le C N^{-1}.\ea\ee
It follows from \eqref{1.1}$_1$ that
\begin{align}\la{oo0}
\frac{d}{dt}\int \n \vp_{N} dx &=\int \n u \cdot\na \vp_{N} dx \notag \\
&\ge -  C N^{-1}\left(\int\n dx\right)^{1/2}\left(\int\n |u|^2dx\right)^{1/2}\ge - \ti C N^{-1}
\end{align}
owing to \eqref{3.1} and \eqref{lvgj1}.
Integrating \eqref{oo0} and choosing $N= N_1\triangleq2N_0+4\tilde CT$, we obtain after using \eqref{1.3} that
\begin{align}\la{p1}
\inf\limits_{0\le t\le T}\int_{B_{N_1}} \n dx&\ge \inf\limits_{0\le t\le T}\int \n \vp_{N_1} dx \notag \\
& \ge \int \n_0 \vp_{N_1} dx-\ti CN_1^{-1}T \notag \\
& \ge \int_{B_{N_0}} \n_0dx-\frac{\ti C T}{2N_0+4\tilde C T} \notag \\
& \ge 1/4.
\end{align}
Hence, it follows from \eqref{p1}, \eqref{3.1} and \eqref{3.i2} that for any $v\in \ti D^{1,2}(\O)$,
\begin{equation}\label{06.2'}
 \| v \bar{x}^{-\eta}\|_{L^\frac{s}{\eta}}\leq C(\eta,s)(\|\n^{1/2}v\|_{L^2}+\|\na v\|_{L^2}),
\end{equation}
where $\eta\in(0,1]$ and $s>2$.
In particular, we deduce from \eqref{06.2'}, \eqref{lvgj1}, and \eqref{lv3.5} that
\begin{equation}\label{06.2}
 \| u \bar{x}^{-\eta}\|_{L^\frac{s}{\eta}}\leq C(\|\n^{1/2}u\|_{L^2}+\|\na u\|_{L^2})\le C.
\end{equation}

Multiplying \eqref{1.1}$_1$ by $\bar{x}^{a}$ and integrating the resulting equation by parts  over $\mathbb{R}^2$  yield  that
\begin{align}
\frac{d}{dt}\int\rho\bar{x}^{a}dx
& \leq C\int\rho| u |\bar{x}^{a-1}\log^{2}(e+|x|^2)dx \notag \\
& \leq C\|\rho\bar{x}^{a-1+\frac{8}{8+a}}\|_{L^{\frac{8+a}{7+a}}}\| u \bar{x}^{-\frac{4}{8+a}}\|_{L^{8+a}} \notag \\
& \leq C\int\rho\bar{x}^{a}dx+C,
\end{align}
which along with  Gronwall's inequality gives \eqref{06.1} and finishes the proof of Lemma \ref{lem03.6}.   \hfill $\Box$

\begin{lemma}\label{lem3.5} There exists a positive constant $C$ depending on $T$ such that
\begin{align}\label{16.1}
\sup_{t\in[0,T]}\|\rho\|_{H^{1}\cap W^{1,q}}
&+\int_{0}^{T}\left(\|\nabla^{2} u \|_{L^2}^2+\|\nabla^{2} u \|_{L^q}^{\frac{q+1}{q}}
+t\|\nabla^{2} u \|_{L^2\cap L^q}^2 \right)dt \notag \\
&+\int_{0}^{T}\left(\|\nabla p \|_{L^2}^2+\|\nabla p \|_{L^q}^\frac{q+1}{q}+t\|\nabla p\|_{L^2\cap L^q}^2\right)dt
\leq C(T).
\end{align}
\end{lemma}

{\it Proof.} First, it follows from
  the mass equation \eqref{1.1}$_1$ that $\na\n$ satisfies  for any $r\ge 2,$
\begin{equation}\label{4.2}
\frac{d}{dt}\|\nabla\rho\|_{L^r}\leq C\|\nabla u \|_{L^\infty}\|\nabla\rho\|_{L^r}.
\end{equation}

Next, one gets from Gagliardo-Nirenberg inequality \eqref{g2}, \eqref{lv3.5}, and \eqref{stokes2} that for $q>2$,
\begin{align}\label{4.3}
\|\na u\|_{L^\infty}
& \leq C(q)\|\nabla u\|_{L^2}^{\frac{q-2}{2(q-1)}}\|\nabla^2u\|_{L^q}^{\frac{q}{2(q-1)}}
\notag \\
& \le C\left(\|\rho\dot{ u }\|_{L^q}^{\frac{q}{2(q-1)}}+\||H||\na H|\|_{L^q}^{\frac{q}{2(q-1)}}\right).
\end{align}

Notice that, it is easy to deduce  from \eqref{3.1}, \eqref{06.2'},   and \eqref{06.1} that for any $\eta\in(0,1]$ and any $s>2$,
\begin{align}\label{lvb01}
\|\rho^\eta v\|_{L^{\frac{s}{\eta}}}
&  \leq C\|\rho^\eta\bar x^{\frac{3\eta a}{4s}}\|_{L^{\frac{4s}{3\eta}}}
\|v\bar  x^{-\frac{3\eta a}{4s}}\|_{L^{\frac{4s}{\eta}}} \notag \\
&  \leq C\|\n\|_{L^\infty}^{\frac{(4s-3)\eta}{4s}}\|\n\bar x^a\|_{L^1}^{\frac{3\eta}{4s}}\left( \|\n^{1/2} v\|_{L^2}+\|\na v\|_{L^2}\right) \notag \\
& \leq C\left( \|\n^{1/2} v\|_{L^2}+\|\na v\|_{L^2}\right),
\end{align}
which together with Gagliardo-Nirenberg inequality yields that
\begin{align}\label{4.6}
\|\rho\dot{ u }\|_{L^q}
& \leq C\|\rho\dot{ u }\|_{L^2}^{\frac{2(q-1)}{q^2-2}}
\|\rho\dot{ u }\|_{L^{q^2}}^{\frac{q(q-2)}{q^2-2}} \notag \\
&  \leq C\|\rho\dot{ u }\|_{L^2}^{\frac{2(q-1)}{q^2-2}}
\left(\|\sqrt{\rho}\dot{ u }\|_{L^2}+\|\nabla\dot{ u }\|_{L^2}\right)^{\frac{q(q-2)}{q^2-2}} \notag \\
&  \leq C \|\sqrt{\rho}\dot{ u }\|_{L^2}
+C\|\sqrt{\rho}\dot{ u }\|_{L^2}^{\frac{2(q-1)}{q^2-2}}
\|\nabla\dot{ u }\|_{L^2}^{\frac{q(q-2)}{q^2-2}}.
\end{align}
This combined  with \eqref{lv3.5} and \eqref{3.13} implies that
\begin{align}\label{4.7}
&  \int_{0}^{T}\left(\|\rho\dot{ u }\|_{L^q}^{\frac{q+1}{q}}
+t\|\rho\dot{ u }\|_{L^q}^{2}\right)dt \notag \\
&  \leq C\int_{0}^{T}\left(\|\sqrt{\rho}\dot{ u }\|_{L^2}^{2}
+t\|\nabla\dot{ u }\|_{L^2}^{2}
+t^{-\frac{q^{3}-q^{2}-2q-1}{q^{3}-q^{2}-2q}}+1\right)dt \notag \\
&  \leq C.
\end{align}
Furthermore, it follows from   \eqref{lv3.5}  and \eqref{lvgj1}  that
\begin{align}\la{lbq-jia8}
& \int_0^T \left(\||H||\na H|\|_{L^q}^{\frac{q+1}{q}}+ t\||H||\na H|\|_{L^q}^2\right)dt \notag \\
& \le C\int_0^T \left(\|\na^2 H\|_{L^{2}}^{1-1/q^2}+ t\|\na^2 H\|_{L^{2}}^{2-2/q}\right)dt \notag \\
& \le C\int_0^T \left(1+t^q+ \|\na^2 H\|_{L^{2}}^2\right)dt \notag \\
& \le C,
\end{align}
where one has used the following estimate
  \be\la{AMSS3.1}\ba \||H||\na H|\|_{L^q}
& \le C(q)\|H\|_{L^2}^{\frac{1}{q}}\|\na H\|_{L^2}\|\na^2 H\|_{L^2}^{\frac{q-1}{q}} \le C\|\na^2 H\|_{L^2}^{\frac{q-1}{q}}
 \ea\ee
 owing to \eqref{g1}, \eqref{lv3.5}, and \eqref{lvgj1}.

The combination of \eqref{4.3}, \eqref{4.7}, and \eqref{lbq-jia8} gives
\be\label{tt}\int_0^T \|\na u\|_{L^\infty}dt\le C. \ee
Thus, applying Gronwall's inequality to \eqref{4.2} yields
\begin{equation}\label{4.8}
\sup_{t\in[0,T]}\|\nabla\rho\|_{L^2\cap L^q}\leq C(T).
\end{equation}

Finally, it is easy to deduce from \eqref{stokes2}, \eqref{4.7}, \eqref{lbq-jia8},  \eqref{lv3.5}, and  \eqref{lv3.5'} that
\begin{align}\label{4.9}
&\int_{0}^{T}\left(\|\nabla^2 u \|_{L^2}^2+\|\nabla^2 u \|_{L^q}^\frac{q+1}{q}+t\|\nabla^2 u \|_{L^2\cap L^q}^2\right)dt \notag \\
&\quad+\int_{0}^{T}\left(\|\nabla p \|_{L^2}^2+\|\nabla p \|_{L^q}^\frac{q+1}{q}+t\|\nabla p\|_{L^2\cap L^q}^2\right)dt\leq C,
\end{align}
which along with \eqref{3.1} and \eqref{4.8} gives \eqref{16.1},  and finishes the proof of Lemma \ref{lem3.5}.  \hfill $\Box$


Now, we give  the following spatial weighted estimate on the gradient of the density, which has been proved in \cite[Lemma 3.6]{lvzh1}. We omit the detailed  proof here for simplicity.

\begin{lemma}\label{lem3.6} There exists a positive constant $C$ depending on $T$ such that
\begin{equation}\label{6.1}
\sup_{t\in[0,T]}\|\rho\bar{x}^{a}\|_{L^{1}\cap H^{1}\cap W^{1,q}}\leq C(T).
\end{equation}
\end{lemma}

Next, by the similar arguments as in \cite{lvhuang,lvshixu}, we shall show the following spatial weighted  estimates of $H$ and $\na H$, which are crucial to  derive  the estimates on the gradients of both $u_t$ and $H_t$.

\begin{lemma}\la{newle}
There exists a positive  constant $C $ depending   on $T$ 
such that
\be
\ba\la{gj10}
\sup_{t\in[0,T]} \| H\bar{x}^{a/2}\|_{L^2}^2 +\int_{0}^{T} \|\na H\bar{x}^{a/2}\|_{L^2}^2dt\le C(T),
\ea
\ee
 and
\be
\ba\la{gj10'}
\sup_{t\in[0,T]} \left(t\|\na H\bar{x}^{a/2}\|_{L^2}^2\right)+\int_{0}^{T} t\|\Delta H\bar{x}^{a/2}\|_{L^2}^2dt\le C(T).
\ea
\ee
\end{lemma}

{\it Proof.} First, multiplying  \eqref{1.1}$_3$ by $H\bar{x}^a$ and integrating the resulting equality  by parts over $\mr^2$   indicate that
\begin{align}\la{mhd8}
& \frac{1}{2}\left(\| H\bar{x}^{a/2}\|_{L^2}^2\right)_t+\nu \|\na H\bar{x}^{a/2}\|_{L^2}^2 \notag \\
& =\frac{\nu}{2}\int |H|^2\Delta\bar{x}^adx+\int H\cdot\na u\cdot H\bar{x}^adx +\frac{1}{2}\int |H|^2u\cdot\na\bar{x}^adx \notag \\
& \triangleq \hat{J}_1+\hat{J}_2+\hat{J}_3,
\end{align}
where
\begin{align*}
\hat{J}_1&\le C\int |H|^2 \bar{x}^a \bar{x}^{-2}\log^{4}(e+|x|^2) dx\leq C\|H\bar{x}^{a/2}\|_{L^2}^2,\\
\hat{J}_2&\le C\int |\na u||H|^2 \bar{x}^a dx\leq C \|\na u\|_{L^2}\|H\bar{x}^{a/2}\|_{L^4}^2\\
&\le C \|H\bar{x}^{a/2}\|_{L^2}\left(\|\na H\bar{x}^{a/2}\|_{L^2}+\|  H\na\bar{x}^{a/2}\|_{L^2}\right)\\
&\le C\|H\bar{x}^{a/2}\|_{L^2}^2+ \frac{\nu}{4}\|\na H\bar{x}^{a/2}\|_{L^2}^2,\\
 \hat{J}_3
 &\leq    C\| H \bar{x}^{a/2}\|_{L^4}\| H \bar{x}^{a/2}\|_{L^2} \|u\bar{x}^{-3/4}\|_{L^{4}}\\
&\leq C\|H\bar{x}^{a/2}\|_{L^2}^2+\frac{\nu}{4}\| \na H \bar{x}^{a/2}\|_{L^2}^2,
\end{align*}
due to \eqref{g1}, \eqref{lv3.5}, and \eqref{06.2}. Then, substituting the above estimates into \eqref{mhd8},   together with Gronwall's inequality, gives \eqref{gj10}.

Now, multiplying  \eqref{1.1}$_3$ by $\Delta H\bar{x}^a$ and integrating the resultant equality  by parts over $\mr^2$  lead  to
\begin{align}\la{AMSS5}
& \frac{1}{2}\left(\|\na H\bar{x}^{a/2}\|_{L^2}^2\right)_t+\nu \|\Delta H\bar{x}^{a/2}\|_{L^2}^2 \notag \\
& \le C\int|\na H| |H| |\na u| |\na\bar{x}^a|dx+C\int|\na H|^2|u| |\na\bar{x}^a|dx+C\int|\na H| |\Delta H| \bar{x}^adx \notag \\
& \quad+C\int |H||\na u||\Delta H|\bar{x}^adx+C\int |\na u||\na H|^2  \bar{x}^adx \notag \\
& \triangleq \sum_{i=1}^5 \tilde{J}_i.
\end{align}
Using Gagliardo-Nirenberg inequality, \eqref{lv3.5}, \eqref{gj10}, and \eqref{06.2}, it holds
\begin{align*}\la{AMSS6}
\tilde{J}_1
\le & C\|H\bar{x}^{a/2}\|_{L^4}^4+C\|\na u\|_{L^4}^4+C\|\na H\bar{x}^{a/2}\|_{L^2}^2\\
\le & C\|H\bar{x}^{a/2}\|_{L^2}^2\left(\|\na H\bar{x}^{a/2}\|_{L^2}^2+\|H\bar{x}^{a/2}\|_{L^2}^2\right)+C\|\na u\|_{L^4}^4+C\|\na H\bar{x}^{a/2}\|_{L^2}^2\\
\le &C+C\|\na^2 u\|_{L^2}^2+C\|\na H\bar{x}^{a/2}\|_{L^2}^2,\\
\tilde{J}_2
\leq  &C\||\na H|^{2-\frac{2}{3a}}\bar{x}^{a-\frac{1}{3}}\|_{L^{\frac{6a}{6a-2}}} \|u\bar{x}^{-\frac{1}{3}}\|_{L^{6a}}\||\na H|^{\frac{2}{3a} }\|_{L^{6a}}\\
\le &C\|\na H \bar{x}^{a/2} \|_{L^2}^\frac{6a-2}{3a}\|\na H \|_{L^4}^\frac{2}{3a}  \leq  C\|\na H \bar{x}^{a/2} \|_{L^2}^2+  C\|\na H \|_{L^4}^2\\
\leq & C\|\na H \bar{x}^{a/2}\|_{L^2}^2+   \frac{\nu}{4}\|\Delta H \bar{x}^{a/2}\|_{L^2}^2,\\
\tilde{J}_3+\tilde{J}_4\le &\frac{\nu}{4}\|\Delta H\bar{x}^{a/2}\|_{L^2}^2+C\|\na H\bar{x}^{a/2}\|_{L^2}^2+C\|H\bar{x}^{a/2}\|_{L^4}^4+C\|\na u\|_{L^4}^4\\
\le &\frac{\nu}{4}\|\Delta H\bar{x}^{a/2}\|_{L^2}^2+C\|\na H\bar{x}^{a/2}\|_{L^2}^2+C\|\na^2 u\|_{L^2}^2+C,\\
\tilde{J}_5\le &C\|\na u\|_{L^\infty} \|\na H\bar{x}^{a/2}\|_{L^2}^2\\
\le &C\left(1+\|\na^2 u\|_{L^q}^{(q+1)/q}\right)\|\na H\bar{x}^{a/2}\|_{L^2}^2.
\end{align*}
Inserting the above estimates into \eqref{AMSS5}  implies that 
\be\la{AMSS10}\ba
 & \left(\|\na H\bar{x}^{a/2}\|_{L^2}^2\right)_t+\nu \|\Delta H\bar{x}^{a/2}\|_{L^2}^2\\
&\quad\le  C\left(1+\|\na^2 u\|_{L^q}^{(q+1)/q}\right) \|\na H\bar{x}^{a/2}\|_{L^2}^2+C(\|\na^2 u\|_{L^2}^2 +1),
\ea\ee
which multiplied by $t$, together with  Gronwall's inequality,
\eqref{gj10},  and \eqref{16.1} yields \eqref{gj10'}.
The proof of Lemma \ref{newle} is finished. \hfill $\Box$

\begin{lemma}\label{lem4.5v}  There exists a positive  constant $C $ depending   on $T$ such that 
\begin{equation}\la{gj13}\ba
\sup_{t\in[0,T]}t\left(\|\n^{1/2}u_t\|^2_{L^2}+\|H_t\|^2_{L^2} +\|\na^2H\|^2_{L^2}\right)+\int_0^T\left(t\|\na u_t\|_{L^2}^2+t\|\na H_t\|_{L^2}^2\right)dt\le C(T).\ea
\end{equation}
\end{lemma}

{\it Proof.}
First, it is easy to deduce from   \eqref{06.2}, \eqref{lvb01}, \eqref{lvgj1}, and \eqref{lv3.5} that
for any $\eta\in(0,1]$ and any $s>2,$
 \be\label{5.d2}\|\n^\eta u \|_{L^{s/\eta}}+ \|u\bar x^{-\eta}\|_{L^{s/\eta}}\le C.\ee

Next, we will prove the following estimate
\be\la{gj12}
\sup_{t\in[0,T]}\left( \|\na u\|_{L^2}^2+ \|\na H\|_{L^2}^2\right)+\int_0^T\left(\|\rho^{1/2}u_t\|_{L^2}^2+ \|H_t\|_{L^2}^2+ \|\Delta H\|_{L^2}^2\right)dt\leq C.
\ee
With \eqref{lv3.5} at hand, we need only to show
\be\la{gj12'} \int_0^T\left(\|\rho^{1/2}u_t\|_{L^2}^2+ \|H_t\|_{L^2}^2\right)dt\leq C.
\ee
Indeed, on the one hand,
\begin{align}\label{lvq1}
\|\rho^{1/2}u_t\|_{L^2}^2
& \le \|\rho^{1/2}\dot u\|_{L^2}^2+\|\rho^{1/2}|u||\na u|\|_{L^2}^2 \notag \\
& \le \|\rho^{1/2}\dot u\|_{L^2}^2+C\|\rho^{1/2} u \|_{L^6}^2\|\na u\|_{L^3}^2 \notag \\
& \le \|\rho^{1/2}\dot u\|_{L^2}^2+C\|\na u\|_{L^2}^{2}+C\|\na^2 u\|_{L^2}^{2}
\end{align}
owing to \eqref{g1} and \eqref{5.d2}. On the other hand,  \eqref{1.1}$_3$ combined with \eqref{g1} and \eqref{lv3.5} leads to
\begin{align}\label{lvq2}
\|H_t\|_{L^2}^2&\le C\|\Delta H\|_{L^2}^2+\||H||\na u|\|_{L^2}^2+\||u||\na H|\|_{L^2}^2 \notag \\
& \le C\|\Delta H\|_{L^2}^2+\|H\|_{L^4}^2\|\na u\|_{L^4}^2+\||u||\na H|\|_{L^2}^2 \notag \\
& \le C\|\Delta H\|_{L^2}^2+C \|\na u\|_{L^2}^2+C\|\na^2 u\|_{L^2}^2+C\| \na H\bar x^{a/2}\|_{L^2}^2,
\end{align}
where in the last inequality one has used
\begin{align}\la{lv4.7}
\||u| |\na H|\|_{L^2}^2
& \leq C\|u \bar{x}^{-a/4}\|_{L^8}^4 \|\na H\|_{L^4}^2
  +C\|\na H \bar{x}^{a/2}\|_{L^2}^2 \notag \\
& \leq \frac{1}{2}\|\na^2  H\|_{L^2}^2+C\|\na H \bar{x}^{a/2}\|_{L^2}^2
\end{align}
due to \eqref{5.d2} and  \eqref{g1}. Hence, \eqref{gj12'} is a direct consequence of  \eqref{lvq1}, \eqref{lvq2}, \eqref{lv3.5}, \eqref{lvgj1}, \eqref{16.1}, and \eqref{gj10}.

Now, differentiating $\eqref{1.1}_2$ with respect to $t$ gives
\begin{align}\la{zb1}
&\n u_{tt}+\n u\cdot \na u_t-\mu\Delta u_t+\na p_t \notag \\ &=-\n_t(u_t+u\cdot\na u)-\n u_t\cdot\na u +\left(H\cdot\na H-\frac{1}{2}\na |H|^2\right)_t.
\end{align}
Multiplying \eqref{zb1} by $u_t$ and integrating the resulting equality by parts over $\O$, we obtain after using $\eqref{1.1}_1$ and $\eqref{1.1}_4$  that
\begin{align}\la{na8}
&\frac{1}{2}\frac{d}{dt} \int \n |u_t|^2dx+\mu\int |\na u_t|^2dx \notag \\
& \leq C\int\rho| u || u_t|(|\nabla u_t|+|\nabla u |^2
+| u ||\nabla^2 u |)dx+ C\int\rho| u |^2|\nabla u ||\nabla u_t|dx \notag \\
& \quad +C\int\rho| u_t|^2|\nabla u|dx-\int H_t\cdot \na u_t\cdot H dx-\int H\cdot \na u_t\cdot H_t dx \notag \\
& \triangleq \sum_{i=1}^{5}\bar{J}_i.
\end{align}

We estimate each term on the right-hand side of \eqref{na8} as follows.

It follows from \eqref{lvb01}, \eqref{5.d2}, \eqref{lv3.5}, \eqref{g1}, and H{\"o}lder's inequality  that
\begin{align}\label{7.5}
\bar{J}_1
& \leq C\|\sqrt{\rho} u \|_{L^6}\|\sqrt{\rho} u_t  \|_{L^2}^{1/2}
\|\sqrt{\rho} u_t  \|_{L^6}^{1/2}
\left(\|\nabla u_t  \|_{L^2}+\|\nabla u \|_{L^4}^2\right) \notag  \\
& \quad +C\|\rho^{1/4} u \|_{L^{12}}^{2}\|\sqrt{\rho} u_t  \|_{L^2}^{1/2}
\|\sqrt{\rho} u_t  \|_{L^6}^{1/2}\|\nabla^2 u \|_{L^2} \notag \\
& \leq C\|\sqrt{\rho} u_t  \|_{L^2}^{1/2}
\left(\|\sqrt{\rho} u_t  \|_{L^2}+\|\nabla u_t  \|_{L^2}\right)^{1/2}
\left( \|\nabla u_t  \|_{L^2}+\|\nabla^2 u \|_{L^2}\right)
\notag \\
& \leq \frac{\mu}{6}\|\nabla u_t \|_{L^2}^{2}+C
\left(1+\|\sqrt{\rho} u_t  \|_{L^2}^{2}
+\|\nabla^2 u \|_{L^2}^2\right).
\end{align}
Next, H{\"o}lder's inequality, \eqref{5.d2}, and \eqref{lvb01} imply
\begin{align}\label{7.6}
\bar{J}_2+\bar{J}_3
& \leq C\|\sqrt{\rho} u \|_{L^8}^{2}\|\nabla u \|_{L^4}\|\nabla u_t  \|_{L^2}
+\|\nabla u \|_{L^2}\|\sqrt{\rho} u_t  \|_{L^6}^{3/2}
\|\sqrt{\rho} u_t \|_{L^2}^{1/2} \notag \\
& \leq \frac{\mu}{6}\|\nabla u_t  \|_{L^2}^{2}
+C \left(1+\|\sqrt{\rho} u_t\|_{L^2}^{2}+\|\nabla^2 u\|_{L^2}^2\right).
\end{align}
For the terms $\bar{J}_4$ and $\bar{J}_5$, using   \eqref{lv3.5} and \eqref{g1},  we obtain   that
\begin{align}\la{na3}
\bar{J}_4+\bar{J}_5
&\le  \frac{\mu}{6} \|\na u_t\|_{L^2}^2+C\||H||H_t|\|_{L^2}^2 \notag \\
&\le  \frac{\mu}{6} \|\na u_t\|_{L^2}^2+C\|H\|_{L^4}^2\|H_t\|_{L^4}^2 \notag \\
&\le  \frac{\mu}{6} \|\na u_t\|_{L^2}^2+C\|H_t\|_{L^2} \|\na H_t\|_{L^2} \notag \\
&\le  \frac{\mu}{6} \|\na u_t\|_{L^2}^2+C\|H_t\|_{L^2}^2+\frac{\mu\nu}{4(C_3+1)} \|\na H_t\|_{L^2}^2,
\end{align}
where $C_3$ is defined in the following \eqref{lv4.13}.
 Submitting \eqref{7.5}--\eqref{na3} into \eqref{na8}  gives
\begin{align}\la{na8-new}
\frac{d}{dt} \|\n^{1/2} u_{t}\|_{L^{2}}^{2}+\mu\|\na u_t\|_{L^2}^2& \le C\left(\|\n^{1/2} u_{t}\|_{L^{2}}^{2}+ \|H_t\|_{L^2}^2\right) \notag \\
&\quad+\frac{\mu\nu}{2(C_3+1)}\|\na H_t\|_{L^2}+C\left( \|\na^2 u\|_{L^2}^2 +1\right).
\end{align}

Next, differentiating $\eqref{1.1}_3$ with respect to $t$ shows
  \be\la{lv4.12}\ba
H_{tt}-H_t\cdot\na u-H\cdot\na u_t+u_t\cdot\na H+u\cdot\na H_t =\nu \Delta H_t.
\ea\ee
 Multiplying \eqref{lv4.12} by $H_t$ and integrating the resulting equality  by parts over $\O,$
 it follows from \eqref{1.1}$_4$, \eqref{g1}, \eqref{lv3.5}, \eqref{06.2'}, and \eqref{gj10} that
\begin{align}\la{lv4.13}
&\frac{1}{2}\frac{d}{dt} \int |H_t|^2dx+\nu \int  |\na H_t|^2dx  \notag \\
&= \int H\cdot\na u_t\cdot H_tdx+\int  H_t\cdot\na u\cdot H_tdx+\int u_t\cdot\na H_t\cdot Hdx  \notag \\
&\leq  C\|H_t\|_{L^4}\|H\|_{L^4}\|\na u_t\|_{L^2}
+C \|H_t\|_{L^4}^2  \|\na u\|_{L^2}  \notag \\
&\quad+C\|H\|_{L^{4}}^{\frac{1}{2a}} \|H \bar x^{a/2}\|_{L^{2}}^{\frac{2a-1}{2a}}\|u_t \bar x^{-\frac{2a-1}{4}}\|_{L^{8a}} \|\na H_t\|_{L^2} \notag \\
&\leq  C\|H_t\|_{L^4}^2+C\|\na u_t\|_{L^2}^2+C\left(\|\n^{1/2} u_{t}\|_{L^{2}} + \|\na u_t\|_{L^2} \right) \|\na H_t\|_{L^2}  \notag \\
&\leq \frac{\nu}{2}\|\na H_t\|_{L^2}^2+C \left(\|H_t\|_{L^2}^2+\|\n^{1/2} u_{t}\|_{L^{2}}^{2}\right)+\frac{C_3}{2}\|\na u_t\|_{L^2}^2.
\end{align}

Now, multiplying \eqref{na8-new}   by $ \mu^{-1}(C_3+1)$  and adding the resulting inequality into   \eqref{lv4.13},  we have
\begin{align}\la{ina8-new}
& \frac{d}{dt}\left( \mu^{-1}(C_3+1)\|\n^{1/2} u_{t}\|_{L^{2}}^{2}+\|H_{t}\|_{L^{2}}^{2}\right)+\|\na u_t\|_{L^2}^2+\frac{\nu}{2} \|\na H_t\|_{L^2}^2 \notag \\
& \le C\left(\|\n^{1/2} u_{t}\|_{L^{2}}^{2}+ \|H_t\|_{L^2}^2\right)+C\left( \|\na^2 u\|_{L^2}^2 +1\right),
\end{align}
which multiplied by $t$, together with  Gronwall's inequality,
\eqref{gj12}, and \eqref{16.1} leads to
\begin{equation}\la{jiagj13}\ba
\sup_{t\in[0,T]}t\left(\|\n^{1/2}u_t\|^2_{L^2}+\|H_t\|^2_{L^2} \right)+\int_0^T\left(t\|\na u_t\|_{L^2}^2+t\|\na H_t\|_{L^2}^2\right)dt\le C(T).\ea
\end{equation}

Finally,  it follows from $\eqref{1.1}_3$, \eqref{g1}, \eqref{lv3.5}, and  \eqref{lv4.7} that
\begin{align}\la{AMSS11}
& \|\na^2 H\|^2_{L^2} \notag \\
& \leq  C\|H_t\|^2_{L^2}+ C\||H| |\na u|\|^2_{L^2}+ C\||u| |\na H|\|_{L^2}^2 \notag \\
& \leq C\|H_t\|^2_{L^2}+C\|H\|^2_{L^4}\|\na u\|_{L^2}\|\na^2 u\|_{L^2}+ \frac{1}{2}\|\na^2  H\|_{L^2}^2+C\|\na H \bar{x}^{a/2}\|_{L^2}^2 \notag \\
& \leq C\|H_t\|^2_{L^2}+C \|\na^2 u\|_{L^2}+\frac{1}{2}\|\na^2  H\|_{L^2}^2+C\|\na H \bar{x}^{a/2}\|_{L^2}^2,
\end{align}
which combined with \eqref{jiagj13}, \eqref{i3.13}, and \eqref{gj10'} indicates \eqref{gj13} and finishes the proof of Lemma \ref{lem4.5v}.
\hfill $\Box$

\section{Proof of Theorem \ref{thm1}}\label{sec4}
With all the a priori estimates in Section \ref{sec3} at hand, we are ready to  prove  Theorem \ref{thm1}.

{\it Proof of Theorem \ref{thm1}.} By Lemma \ref{lem21}, there exists a $T_{*}>0$ such that the problem \eqref{1.1}-\eqref{1.3'} has a unique strong solution $(\rho, u ,p, H)$ on $\mathbb{R}^2\times(0,T_{*}]$. Now, we will extend the local solution to all time.


Set
\begin{equation}\label{20.1}
T^{*}=\sup \left\{T~|~(\rho,u,p, H)\ \text{is a strong solution on}\ \mathbb{R}^{2}\times(0,T]\right\}.
\end{equation}
First, for any $0<\tau<T_*<T\leq T^{*}$ with $T$ finite, one deduces from \eqref{lv3.5}, \eqref{lvgj1},  \eqref{i3.13}, and \eqref{gj13}  that for any $q\geq2$,
\begin{equation}\label{20.2}
\nabla u, \na H, H \in C([\tau,T];L^2\cap L^q),
\end{equation}
where one has used the standard embedding
\begin{equation*}
L^{\infty}(\tau,T;H^1)\cap H^{1}(\tau,T;H^{-1})\hookrightarrow C(\tau,T;L^q)\ \ \text{for any}\ \ q\in[2,\infty).
\end{equation*}
Moreover, it follows from  \eqref{16.1}, \eqref{6.1}, and \cite[Lemma 2.3]{L1996} that
\begin{equation}\label{20.3}
\rho\in C([0,T];L^1\cap H^1\cap W^{1,q}).
\end{equation}

Finally, we claim that
\begin{equation}\label{20.4}
T^{*}=\infty.
\end{equation}
Otherwise, if $T^{*}<\infty$,      it follows from \eqref{20.2},  \eqref{20.3},  \eqref{lv3.5}, \eqref{lvgj1},     \eqref{6.1},   and \eqref{gj10}  that
$$(\n, u, H)(x,T^*)=\lim_{t\rightarrow T^*}(\n, u,H)(x,t)$$
satisfies the initial conditions \eqref{2.2} at $t=T^*$. Thus, taking $(\n, u,H)(x,T^*)$ as the initial data, Lemma 2.1 implies that one could extend the local strong solutions beyond $T^*$. This  contradicts the assumption  of $T^*$ in \eqref{20.1}.  The proof of Theorem \ref{thm1} is completed.  \hfill $\Box$


\begin{thebibliography}{10}

 \bibitem{abidi1} H. Abidi and T. Hmidi,
 \emph{R\'esultats d\'existence dans des espaces critiques pour le syst\'eme de la MHD inhomog\'ene}, Ann. Math. Blaise Pascal, \textbf{14}(2007), 103--148.


\bibitem{abidi2} H. Abidi and M. Paicu,
\emph{Global existence for the MHD system in critical spaces},  Proc. Roy. Soc. Edinburgh  Sect. A, \textbf{138}(2008), 447--476.


\bibitem{AK1973}
\newblock S. A.~Antontesv and A. V.~Kazhikov,
\newblock  Mathematical Study of Flows of Nonhomogeneous Fluids ,
\newblock Lecture notes, Novosibirsk State University, Novosibirsk, U.S.S.R., 1973(in Russian).


\bibitem{AKM1990}
\newblock S. A.~Antontesv, A. V.~Kazhikov and V. N. Monakhov,
\newblock Boundary Value Problems in Mechanics of Nonhomogeneous Fluids ,
\newblock North-Holland, Amsterdam, 1990.


\bibitem{tan} Q. Chen, Z. Tan and Y. J. Wang,
\emph{Strong solutions to the incompressible magnetohydrodynamic equations}, Math. Methods Appl. Sci., \textbf{34(1)}(2011), 94--107.


\bibitem{CK2003}
\newblock H. J.~Choe and H.~Kim,
\newblock \emph{Strong solutions of the Navier-Stokes equations for nonhomogeneous incompressible fluids},
\newblock Comm. Partial Differential Equations, \textbf{28} (2003), 1183--1201.


\bibitem{CLMS1993}
\newblock R.~Coifman, P. L.~Lions, Y.~Meyer and S.~Semmes,
\newblock \emph{Compensated compactness and Hardy spaces},
\newblock J. Math. Pures Appl., \textbf{72} (1993), 247--286.


\bibitem{david} P. A. Davidson,
An Introduction to Magnetohydrodynamics, Cambridge University
Press, Cambridge, 2001.


\bibitem{D1997}
\newblock B.~Desjardins,
\newblock \emph{Regularity results for two-dimensional flows of multiphase viscous fluids},
\newblock Arch. Rational Mech. Anal., \textbf{137} (1997), 135--158.


\bibitem{leb2} B. Desjardins and C. Le Bris,
\emph{Remarks on a nonhomogeneous model of magnetohydrodynamics},  Differential Integral Equations, \textbf{11(3)}(1998), 377-394.


\bibitem{dlions}
G. Duraut and J. L. Lions,
 \emph{In\'equations en thermo\'elasticit\'e et magn\'etohydrodynamique},
  Arch. Rational Mech. Anal., \textbf{46} (1972), 241-279.


\bibitem{F2004}
\newblock E.~Feireisl,
\newblock  Dynamics of Viscous Compressible Fluids ,
\newblock Oxford University Press, Oxford, 2004.


\bibitem{leb1} J. F. Gerbeau and C. Le Bris,
\emph{Existence of solution for a density-dependant magnetohydrodynamic equation}, Adv. Differential Equations,
\textbf{2(3)}(1997), 427--452.


 \bibitem{H19951}
\newblock D.~Hoff,
\newblock \emph{Global solutions of the Navier-Stokes equations for multidimensional compressible
flow with discontinuous initial data},
 \newblock J. Differential Equations, \textbf{120} (1995), 215--254.


\bibitem{HL20130}
\newblock X. D.~Huang and J.~Li,
\newblock \emph{Global classical and weak solutions to the three-dimensional full compressible Navier-Stokes system with vacuum and large oscillations},
\newblock http://arxiv.org/abs/1107.4655 [math-ph], 2011.


\bibitem{HL20131}
\newblock X. D.~Huang and J.~Li,
\newblock \emph{Existence and blowup behavior of global strong solutions to the two-dimensional baratropic compressible Navier-Stokes system with vacuum and large initial data},
\newblock http://arxiv.org/abs/1205.5342 [math.AP], 2012.


\bibitem{hwjde1}
\newblock X. D.~Huang and Y.~Wang,
\newblock \emph{Global strong solution to the 2D nonhomogeneous incompressible MHD system},
\newblock J. Differential Equations, \textbf{254}(2014), 511--527.


\bibitem{HW2014}
\newblock X. D.~Huang and Y.~Wang,
\newblock \emph{Global strong solution with vacuum to the two-dimensional density-dependent Navier-Stokes system},
\newblock SIAM J. Math. Appl., \textbf{46} (2014), 1771--1788.


\bibitem{HW2015}
\newblock X. D.~Huang and Y.~Wang,
\newblock \emph{Global  strong  solution  of  3D  inhomogeneous Navier-Stokes  equations with  density-dependent viscosity},
\newblock J. Differential Equations, \textbf{259}(2015), 1606--1627.


\bibitem{HLX2012}
\newblock X. D.~Huang, J.~Li and Z. P.~Xin,
\newblock \emph{Global well-posedness of classical solutions with large oscillations and vacuum to the three-dimensional isentropic compressible Navier-Stokes equations},
\newblock Comm. Pure Appl. Math., \textbf{65} (2012), 549--585.


\bibitem{K1986}
\newblock T.~Kato,
\newblock \emph{Strong $L^p$-solutions of the Navier-Stokes equation in $\mathbb{R}^m$, with applications to weak solutions}.,
\newblock Math. Z. \textbf{187}(4)(1984), 471-480.


\bibitem{K1974}
\newblock A. V.~Kazhikov,
\newblock \emph{Resolution of boundary value problems for nonhomogeneous viscous fluids},
\newblock Dokl. Akad. Nauk., \textbf{216} (1974), 1008--1010.


\bibitem{LL2014}
\newblock J.~Li and Z. L.~Liang,
\newblock \emph{On Local classical solutions to the Cauchy problem of the two-dimensional
barotropic compressible Navier-Stokes equations with vacuum},
\newblock J. Math. Pures Appl., \textbf{102} (2014), 640--671.


\bibitem{LX2014}
\newblock J.~Li and Z. P. Xin,
\newblock \emph{Global well-posedness and large time asymptotic behavior of classical solutions to the compressible Navier-Stokes equations with vacuum},
\newblock http://arxiv.org/abs/1310.1673 [math.AP], 2013.


\bibitem{L2015}
\newblock Z. L.~Liang,
\newblock \emph{Local strong solution and blow-up criterion for the 2D nonhomogeneous incompressible fluids},
\newblock J. Differential Equations, \textbf{7}(2015),  2633--2654.


\bibitem{L1996}
\newblock P. L.~Lions,
\newblock  Mathematical Topics in Fluid Mechanics , vol. I: Incompressible Models,
\newblock Oxford University Press, Oxford, 1996.


\bibitem{lvzh1}
\newblock B. Q. L\"u, X. D. Shi and X. Zhong,
\newblock \emph{Global existence and large time asymptotic behavior of strong solutions to the Cauchy problem of 2D density-dependent  Navier-Stokes equations with vacuum},
http://arxiv.org/abs/1506.03143 [math.AP], 2015.


\bibitem{lvzh}
\newblock B. Q. L\"u, Z. H. Xu and X. Zhong,
\newblock \emph{On local strong solutions to the Cauchy problem of the two-dimensional density-dependent magnetohydrodynamic equations with vacuum},
http://arxiv.org/abs/1506.02156 [math.AP], 2015.

\bibitem{lvhuang}
\newblock B. Q. L\"u and B. Huang,
\newblock \emph{On strong solutions to the Cauchy problem of the two-dimensional compressible magnetohydrodynamic equations with vacuum },
\newblock Nonlinearity, \textbf{28}(2015), 509--530.


\bibitem{lvshixu}
\newblock B. Q. L\"u, X. D. Shi and X. Y. Xu,
\newblock \emph{Global  well-posedness  and large time asymptotic behavior of strong solutions to the  compressible magnetohydrodynamic equations with vacuum}, \newblock http://arxiv.org/abs/1402.4851 [math.AP], 2014.


\bibitem{nir} L. Nirenberg,
\emph{On elliptic partial differential equations}.  Ann. Scuola Norm. Sup. Pisa, (3){\bf 13}(1959), 115--162.


\bibitem{teman} M. Sermange and R. Temam,
\emph{Some mathematical questions related to the MHD equations}. Comm.  Pure Appl.  Math.,   \textbf{36}(1983), 635--664.


\bibitem{S1990}
\newblock J.~Simon,
\newblock \emph{Nonhomogeneous viscous incompressible fluids: existence of velocity, density, and pressure},
\newblock SIAM J. Math. Anal., \textbf{21} (1990), 1093--1117.


\bibitem{S1993}
\newblock Elias M.~Stein,
\newblock Harmonic Analysis: Real-variable Methods, Orthogonality, and Oscillatory Integrals,
\newblock Princeton University Press, Princeton, NJ, 1993.


\bibitem{T2001}
\newblock R.~Temam,
\newblock Navier-Stokes Equations: Theory and Numerical Analysis.
Reprint of the 1984 edition,
\newblock AMS Chelsea Publishing, Providence, RI, 2001.


\bibitem{zhang}
\newblock J. W.~Zhang,
\newblock \emph{Global well-posedness  for  the  incompressible  Navier-Stokes  equations with  density-dependent  viscosity  coefficient},
\newblock J. Differential Equations, \textbf{259(5)}(2015), 1722--1742.

\end{thebibliography}
\end{document}